\documentclass[11pt]{amsart}
\usepackage{amssymb,amsmath}
\usepackage{mathrsfs}

\usepackage{color}

\numberwithin{equation}{section}

\def\H{H^\infty}

\newcommand{\calA}{\mathcal{A}}

\newcommand{\mC}{\mathbb{C}}
\newcommand{\mD}{\mathbb{D}}

\newcommand{\mN}{\mathbb{N}}

\newtheorem{theorem}{Theorem}[section]
\newtheorem{lemma}[theorem]{Lemma}
\newtheorem{corollary}[theorem]{Corollary}
\newtheorem{proposition}[theorem]{Proposition}
\newtheorem{conjecture}[theorem]{Conjecture}

\theoremstyle{definition}

\theoremstyle{definition}
\newtheorem{definition}[theorem]{Definition}

\begin{document}



\title[A Corona theorem]{The Corona theorem
 and stable rank for the algebra $\mC +BH^\infty$}


\author{Raymond Mortini}
\address{Mortini:
 Universit\'e Paul Verlaine - Metz,
 D\'epartement de Math\'ematiques,
 Ile du Saulcy,
 F-57045 METZ
 France.}
\email{mortini@univ-metz.fr}

\author[Amol Sasane]{Amol Sasane$^*$}
\thanks{$^*$ Research supported by the Nuffield Grant NAL/32420.}
\address{Sasane:
Department of Mathematics,
London School of Economics,
Houghton Street,
London WC2A 2AE,
United Kingdom.}
\email{A.J.Sasane@lse.ac.uk}

\author[Brett D. Wick]{Brett D. Wick$^\dagger$}
\thanks{$\dagger$ Research supported in part by the National Science Foundation under DMS Grant \# 0752703.}
\address{Wick:
Department of Mathematics,
University of South Carolina,
LeConte College,
1523 Greene Street,
Columbia, SC 29208, USA.
}
\address{
\hspace*{3em}Present:
Department of Mathematics,
The Royal Institute of Technology (KTH),
S -- 100 44 Stockholm, Sweden.}
\email{wick@math.sc.edu}

\subjclass[2000]{Primary  46J15 ; Secondary 32A38, 19B10, 46J30, 30H05}


\keywords{Subalgebras of $H^\infty$, Corona theorem, Bezout equation, Bass stable rank,
Ideals, Blaschke products}

\begin{abstract}
  Let $B$ be a Blaschke product. We prove in several different ways
  the corona theorem for the algebra $H^\infty_B:=\mC+BH^\infty$.
  That is, we show the equivalence of the classical {\em corona
    condition} on data $f_1, \dots, f_n \in H^\infty_B$:
  \[
  \forall z \in \mD, \;\; \sum_{k=1}^{n} |f_k(z)| \geq \delta >0,
  \]
  and the {\em solvability of the Bezout equation} for $g_1, \dots,
  g_n \in H^\infty_B$:
  \[
  \forall z\in \mD, \;\; \sum_{k=1}^n g_k (z)f_k(z)=1.
  \]
  Estimates on solutions to the Bezout equation are also obtained.  We
  also show that the Bass stable rank of $H^\infty_B$ is $1$.  Let
  $A(\mD)_B$ be the subalgebra of all elements from $H^\infty_B$ having a
  continuous extension to the closed unit disk $\overline{\mD}$.
  Analogous results are obtained also for $A(\mD)_B$.
\end{abstract}
 \maketitle

\section{Introduction}

Let $\mD:=\{z\in \mC\;:\; |z|<1\}$.  We denote the Hardy algebra of
all  bounded holomorphic functions on $\mD$ (with
pointwise  operations and the supremum norm) by $H^\infty$, while the
disk algebra is the set of all functions in $H^\infty$ having a
continuous extension to $\overline{\mD}$ will be denoted by $A(\mD)$.

It is well known that the unit disk $\mD$ is dense in the maximal
ideal space of $H^\infty$.  This can be rephrased in terms of function
theoretic conditions to say that for any $f_1,\ldots, f_n\in H^\infty$
such that
$$
\forall z\in\mD, \;\;\sum_{k=1}^n|f_k(z)| \geq\delta >0,
$$
there exists $g_1,\ldots, g_n\in H^\infty$ such that
$$
\forall z\in\mD,\;\;\sum_{k=1}^n g_k(z) f_k(z)=1.
$$
It is immediate to see that if $1$ is in the ideal generated by the
functions $f_1,\ldots, f_n$, then the condition
$$
\exists \delta>0 \textrm{ such that }\forall z\in\mD,
\;\;\sum_{k=1}^n|f_k(z)|\geq\delta>0
$$
is necessary.  The fact that this condition is sufficient is the
famous Carleson Corona Theorem \cite{Car62}.  In fact, Carleson proved
that there are estimates on the solutions $g_j$ in terms of the
parameters $\delta$ and $n$, whenever the corona data $f_j$ are in the
unit ball of $H^\infty$.

The goal of this paper is to prove analogous results for certain
subalgebras of $H^\infty$ that arise from constraints on the
functions, and their derivatives at some points in the disk.  With
this in mind, we define the space $H^\infty_B$.

For $a\in\mD$, let
\[
\varphi_a(z):=\frac{a-z}{1-\overline{a}z}
\]
denote the M\"obius transform of the unit disk $\mD$ onto itself.
Then we define the {\em Blaschke product} $B$ with zeros $a_k$ of
multiplicity $m_k$ by
\[
B:=\prod_{k\geq 1} \biggl( \frac{|a_k|}{a_k}\varphi_{a_k} \biggr)^{m_k}
 \textrm{ where }
\sum_{k\geq 1} m_k(1-|a_k|) <\infty.
\]
Note that  the factor
$\frac{|a_k|}{a_k}$ is chosen for convergence of the infinite product
(if $a_k=0$ then this factor is $-1$).

\begin{definition}
  Let $a_k$, $k\geq 1$, denote the zeros in $\mD$ of a Blaschke
  product $B$ with multiplicity $m_k$.  We denote by $H^\infty_B$ the
  set of all those functions in $H^\infty$ that satisfy
  \begin{itemize}
  \item[(A1)] For all $j,k$: $f(a_j)=f(a_k)$;
  \item[(A2)] For all $k$ and all $1\leq m\leq m_k-1$: $f^{(m)}(a_k)=0$.
  \end{itemize}
  Similarly, $A(\mD)_B$ is the set of all those functions in $A(\mD)$ satisfying
  (A1) and (A2).
\end{definition}
We note that $A(\mD)_B$ reduces to $\mC$ if the Blaschke product
$B$ has its zeros accumulating at a set of positive Lebesgue measure.

An equivalent and useful definition of the algebra $\mC+B\H$ will be
the following:
$$
H^\infty_B:=\mathbb{C}+BH^\infty
$$
We will use the definitions interchangeably.

One notes that $H^\infty_B$ and $A(\mD)_B$ are Banach algebras with the
usual pointwise operations and the supremum norm:
\[
\|f\|_{\infty}:=\sup_{z\in \mD} |f(z)|.
\]
Let us point out that for infinite Blaschke products, $\mC+ BA(\mD)$ is not
an algebra and that $A(\mD)_B$ is strictly contained in $\mC+ BA(\mD)$. For
finite Blaschke products both sets coincide.

The algebra $H^\infty_B$ has been studied earlier in \cite{Rag}. The
case when the Blaschke product has a finite number of simple zeros was
considered in \cite{Sol}, and the simplest case of one zero at $0$
with multiplicity $2$ was discussed in \cite{DavPauRagSin}.  These
papers considered Nevanlinna-Pick interpolation problems for these
algebras.

These algebras also arise naturally as the restriction of functions in
$H^\infty$ and the ball algebra on the bidisc to distinguished
varieties.  In this context, they have been extensively studied by
Agler and McCarthy, see \cite{AM1} and \cite{AM2}.  Additionally, in
the case of a finite Blaschke product, these algebras give examples of
finite co-dimensional subalgebras of $H^\infty$, which were first
studied by Gamelin, see \cite{Gamelin}.

\subsection{Main Results}

In this article, we prove in four different ways a corona theorem for
$H^\infty_B$. The first version reads as follows:

\begin{theorem}
\label{corona_H_infty_B}
The set of multiplicative linear functionals corresponding to
evaluations at points of the unit disk is dense in the maximal ideal
space of $H^\infty_B$.
\end{theorem}

Note that one proof follows immediately from \cite[p.12]{gam} and the
classical Carleson Corona Theorem, which says that $\mD$ is dense in
the maximal ideal space of $H^\infty$ \cite{Car62}.)

We also have the following:

\begin{theorem}
\label{coronaAB}
If $A(\mD)_B\not=\mC$, then the maximal ideal space $M$ of $A(\mD)_B$ is the
set of multiplicative linear functionals corresponding to evaluations
at points of $\overline{\mD}$.  $M$ can be identified with the
quotient space $\overline{\mD}|_{\sim}$, where $z\sim w$ if and only
if $z,w\in \overline{\{a_k: k\in \mN\}}$.
\end{theorem}

We point out here that because the algebras $H^\infty_B$ and $A(\mD)_B$ by
definition have points where the functions agree, the unit disk does
not embed homeomorphically into the spectrum of these algebras. This
is because the points where the functions agree give rise to the same
linear functionals.

By the standard Gelfand theory of Banach algebras these results imply:
\goodbreak

\begin{corollary}
\label{corollary_main}
Let $\calA= H^\infty_B$ or $\calA=A(\mD)_B$. Then the following are
equivalent:
\begin{enumerate}
\item $f_1, \dots , f_n \in \calA$ and there exists a $ \delta>0$, such
  that
\[
\forall z \in \mD, \;\; \sum_{k=1}^{n} |f_k(z)| \geq \delta >0 \quad
\textrm{\em (Corona condition)}.
\]
\item There exist $g_1, \dots, g_n \in \calA$, such that
\[
\forall z\in \mD, \;\; \sum_{k=1}^n g_k (z)f_k(z)=1\quad
\textrm{\em (Bezout equation)}.
\]
\end{enumerate}
\end{corollary}

The well known Carleson Corona Theorem for $H^\infty$ is stronger in
that it is possible to find estimates for the solution to the Bezout
equation in terms of the parameter $\delta$.  It is also possible to
do this in the case of the algebras $H^\infty_B$. We first present a
``matricial'' proof interesting in its own right when $B$ is a
\textit{finite} Blaschke product.

\begin{theorem}
\label{theorem_estimates}
Let $B$ be a finite Blaschke product with zeros $a_k$ of multiplicity
$m_k$, $k=1,\ldots, N$.  Let $f_1, \dots , f_n \in H^\infty_B$ be such
that
\[
\forall z \in \mD, \;\; 1\geq \sum_{k=1}^{n} |f_k(z)| \geq \delta >0.
\]
Then there exist $g_1, \dots, g_n \in H^\infty_B$ such that
\[
\label{Bezout_equation}
\forall z\in \mD, \; \sum_{k=1}^n g_k (z)f_k(z)=1 \;\; \textrm{ and }\;\;
 \forall k \in \{1, \dots,
n\},\; \|g_{k}\|_{\infty} \leq C(n, \delta,a_k,m_k).
\]
\end{theorem}

Using a short argument relying on an application of Carleson's Corona
Theorem, we are able to show that the above theorem holds more
generally for algebras generated by ideals in $H^\infty$. This allows
us to remove the dependence upon the points $a_k$ corresponding to the
zeros of the Blaschke factors.

Namely, let $\mathbb{I}$ be any ideal in $H^\infty$, and let
\[
H^\infty_{\mathbb{I}}:=
\{c+\varphi\;|\; c\in \mC \textrm{ and } \varphi\in \mathbb{I}\}.
\]
Then $H^\infty_{\mathbb{I}}$ is a subalgebra of $H^\infty$.
We then establish the following:

\begin{theorem}
\label{theorem_estimates2}
Let $\mathbb{I}$ be a proper ideal in $H^\infty$ and let $f_1, \dots ,
f_n \in H^\infty_{\mathbb{I}}$ be such that
\begin{equation}
\label{CC_I}
\forall z \in \mD, \;\; 1\geq \sum_{k=1}^{n} |f_k(z)| \geq \delta >0.
\end{equation}
Then there exist $g_1, \dots, g_n \in H^\infty_{\mathbb{I}}$ such that
\[
\forall z\in \mD, \; \sum_{k=1}^n g_k (z)f_k(z)=1 \;\; \textrm{ and
}\;\; \forall k \in \{1, \dots, n\},\; \|g_{k}\|_{\infty} \leq
C(n,\delta).
\]
\end{theorem}

We prove Theorem \ref{theorem_estimates} and Theorem
\ref{theorem_estimates2} in Section \ref{section_estimates} of this
article.  In Section \ref{stableranks}, we show that the Bass stable
rank of $A(\mD)_B$ and $H^\infty_B$ is equal to $1$.

\section{Corona theorem for $H^\infty_B$ and $A(\mD)_B$}

\subsection{The Case of $H^\infty_B$}
Let
\[
H_{B0}^{\infty}:=\{f\in H^\infty_B \; : \;f(a_1)=0\}.
\]
Note that $H_{B0}^{\infty} \subset H_B^{\infty}$, and
$H_{B0}^{\infty}$ is a Banach algebra without identity. It is
straightforward to check that $H_{B0}^{\infty}$ is nothing but the
closed ideal $B\H$ generated by the Blaschke product $B$.

Evaluation at a point $\lambda \in \mD$ is a multiplicative linear
functional on $H_{B0}^{\infty}$ (the trivial one if $\lambda$ belongs
to the zero set $\Lambda$ of $B$).

We prove below that in the weak-* topology on the dual space of
$H_{B0}^{\infty}$, the set of these point evaluations is dense in the
set of all multiplicative linear functionals on $H_{B0}^{\infty}$.

\begin{lemma}
  The set of multiplicative linear functionals consisting of evaluations
  at points of $\mD\setminus \Lambda$ is dense in the set of all
  non-zero multiplicative linear functionals on $H_{B0}^{\infty}$.
\end{lemma}
\begin{proof} Let $l$ be a non-zero multiplicative linear functionals
on $H_{B0}^{\infty}$.  Since $l$ is non-zero, there exists a
function $p_0 \in H_{B0}^{\infty}$ with $l(p_0)\neq 0$. If $f\in
H^\infty$, then define
\[
L(f):=\frac{l(fp_0)}{l(p_0)}.
\]
Note that $(fp_0)(a_k)=f(a_k) 0=0$ and for $1\leq m \leq m_k-1$,
\[
(fp_0)^{(m)}(a_k)=
f^{(m)}(a_k) \underbrace{p_0 (a_k)}_{=0}
+
\sum_{j=1}^{m} \binom{m}{j} f^{(m-j)}(a_k) \underbrace{p_0^{(j)} (a_k)}_{=0}
=0.
\]
So $fp_0\in H_{B0}^{\infty}$, and $L$ is well-defined.  Clearly $L$ is
a linear transformation.  For $f\in H^\infty$, we have
\[
|L(f)|= \frac{|l(f p_0)|}{|l(p_0)|}
\leq \frac{\|fp_0\|_{\infty}}{|l(p_0)|}
\leq
 \bigg( \frac{\|p_0\|_{\infty}}{|l(p_0)|} \bigg) \|f\|_{\infty},
\]
and so $L$ is continuous on $H^\infty$. It is also multiplicative,
since if $f,g\in H^\infty$, then $ l(fg p_0) l(p_0)= l(fg p_0 p_0)=
l((fp_0)(gp_0))=l(fp_0) l(gp_0)$, and so dividing by $(l(p_0))^2$
($\neq 0$), we obtain
\[
L(fg)= \frac{l(fg p_0)}{l(p_0)}=\frac{l(f
  p_0)}{l(p_0)}\frac{l( g p_0)}{l(p_0)}=L(f)L(g).
\]
Thus $L$ defines a non-zero multiplicative linear functional on
$H^\infty$. By the Carleson Corona Theorem, there exists a net
$(\alpha_j)_{j\in J}$ of point evaluations in $\mD$ that
converges to $L$ in the Gelfand topology of the maximal
ideal space of $H^\infty$. We also observe
that $l$ is the restriction of $L$ to $H_{B0}^{\infty}$:
\[
\forall f\in H_{B0}^{\infty}, \;\;L(f)= \frac{l( f
  p_0)}{l(p_0)}=\frac{l(f)l( p_0)}{l(p_0)}=l(f).
\]
The restriction of each element in the net $(\alpha_j)_{j\in J}$ to
$H_{B0}^{\infty}$ gives a net (of point evaluations in $\mD$) that
converges to $l$ in the weak-* topology of $H_{B0}^{\infty}$.
\end{proof}

\begin{proof}[Proof of Theorem \ref{corona_H_infty_B}]
  Let $L$ be a non-zero multiplicative linear functional on
  $H_B^{\infty}$. Let $l:=L|_{H_{B0}^{\infty}}$. Then $l$ is a
  multiplicative linear functional on $H_{B0}^{\infty}$. If $f\in
  H_B^{\infty}$, then $f-f(a_1) \in H_{B0}^{\infty}$, and so
  $L(f)=f(a_1)+l(f-f(a_1))$. Now, if $l$ is identically zero, we have
  $L(f)=f(a_1)$, and so $L$ is point evaluation (at $a_1$).  On the
  other hand, if $l$ is non-zero, then by the previous lemma, there
  exists a net $(\alpha_j)_{j\in J}$ of point evaluations in $\mD$
  that converges to $l$ in the weak-* topology of $H_{B0}^{\infty}$.
  Thus for all $f\in H^\infty_B$,
\begin{eqnarray*}
L(f)
&=& f(a_1)+l(f-f(a_1))
= f(a_1)+\big(\lim_{j\in J} \alpha_j\big)(f-f(a_1))\\
&=& f(a_1)+\lim_{j\in J} \alpha_j \big(f-f(a_1)\big)
= f(a_1)+\lim_{j\in J} ( f(\alpha_j)-f(a_1))\\
&=&\lim_{j\in J} f(\alpha_j)
= \lim_{j\in J} \alpha_j (f).
\end{eqnarray*}
Thus $L=\displaystyle \lim_{j\in J} \alpha_j$, and this completes the proof.
\end{proof}

\subsection{The Case of $A(\mD)_B$}

The proof for $A(\mD)_B$ is basically the same, however is technically
easier since the maximal ideal space of the disk algebra is
$\overline{\mD}$, so the language of nets doesn't need to be employed.

Let
\[
A(\mD)_{B0}:=\{f\in A(\mD)_B \; : \;f(a_1)=0\}.
\]
Note that $A(\mD)_{B0} \subset A(\mD)_B$, and $A(\mD)_{B0}$ is a Banach algebra
without identity.  Moreover $A(\mD)_{B0}= B\H\cap A(\mD)$. In order to avoid
trivialities, we assume that the zeros of $B$ do not cluster at a set
of positive Lebesgue measure, so that $A(\mD)_{B0}$ does not collapse to
$\{0\}$.

\begin{lemma}
  Each non-zero multiplicative linear functional $l:A(\mD)_{B0}
  \rightarrow\mC$ is a point evaluation at some $\lambda\in
  \overline\mD\setminus \overline\Lambda$.
\end{lemma}
\begin{proof}
  Since $l$ is non-zero, there exists a function $p_0 \in A(\mD)_{B0}$ with
  $l(p_0)\neq 0$. For $f\in A(\mD)$, define $L(f)=l(fp_0)/l(p_0)$.  Then
  $L$ defines a non-zero multiplicative linear functional on $A(\mD)$, and
  hence is point evaluation at some $\lambda \in \overline{\mD}$.  As
  $l$ is the restriction of $L$ to $A(\mD)_{B0}$, $l$ is point evaluation
  as well.
\end{proof}

\begin{proof} [Proof of Theorem \ref{coronaAB}]
  Let $L$ be a non-zero multiplicative linear functional on $A(\mD)_B$.
  Define $l:=L|_{A(\mD)_{B0}}$. Then $l$ is a multiplicative linear
  functional on $A(\mD)_{B0}$. If $f\in A(\mD)_B$, then
  $L(f)=f(a_1)+l(f-f(a_1))$. Now, if $l$ is identically zero, we have
  $L(f)=f(a_1)$, and so $L$ is point evaluation (at $a_1$).  On the
  other hand, if $l$ is non-zero, then by the previous lemma, $l$ is
  point evaluation, at say, $\lambda \in \overline{\mD}$.  Thus for
  all $f\in A(\mD)_B$, $L(f)=f(a_1)+l(f-f(a_1))=f(\lambda)$.
\end{proof}

\section{Corona theorem for $H^\infty_B$ with estimates}
\label{section_estimates}

In this section, we show that when the number of zeros of the Blaschke
product $B$ is finite, we can get estimates on the size of a solution
to the Bezout equation in Corollary \ref{corollary_main} for
$H^\infty_B$.  We also give a different proof that is valid for subalgebras
of $H^\infty$ that are generated by ideals, which uses only Carleson's Corona Theorem.

We will use the following simple lemma. Here we use the notation
$M^{\top}$ for the {\em transpose} of a matrix $M\in \mC^{m\times n}$,
i.e., $[M^{\top}]_{ij}=[M]_{ji}$. For a (column)-vector $x\in \mC^n$,
let $||x||$ be the usual Euclidean norm of $x$; the norm $||A||$ of a
matrix $A=(a_{ij})\in \mC^{n\times n}$ is given by
$||A||=(\sum_{i,j}|a_{ij}|^2)^{1/2}$.

Also, in this section we will sometimes use the letter $C$ to denote
different constants in the same proof and we will keep track of the
parameters that determine the constant.

\begin{lemma}
\label{lemma_LA}
Let $x, y \in \mC^{n}$ such that $\|x\| \geq \delta>0$
and $x^{\top} y=0$. Then there exists
$A\in \mC^{n\times n}$ such that $
A=-A^{\top}$, $Ax=y$ and $\|A\|\leq C(n,\delta) \|y\|$.
\end{lemma}
\begin{proof} Since $x\neq 0$, there exists $k\in \{1, \dots, n\}$
  such that $x_{k} \neq 0$ and $|x_{k}|= \max\{ |x_1|, \dots,
  |x_n|\}\geq C(n)\|x\| \geq C(n) \delta$; for example we may take
  $C(n)=1/\sqrt n$.  Define
\[
A:=\left[ \begin{array}{ccc|c|ccc}
     0        & \dots  & 0            & \frac{y_1}{x_k}     & 0            & \dots  & 0        \\
     \vdots   & \ddots & \vdots       & \vdots      & \vdots       & \ddots & \vdots   \\
     0 & \dots & 0 & \frac{y_{k-1}}{x_k} & 0 & \dots & 0 \\ \hline -
     \frac{y_1}{x_k} & \dots & - \frac{y_{k-1}}{x_k} & 0 & -
     \frac{y_{k+1}}{x_k} & \dots & - \frac{y_n}{x_k} \\ \hline
     0        & \dots  & 0            &  \frac{y_{k+1}}{x_k} & 0            & \dots  & 0        \\
     \vdots   & \ddots & \vdots       & \vdots      & \vdots       & \ddots & \vdots   \\
     0 & \dots & 0 & \frac{y_{n}}{x_k} & 0 & \dots & 0
 \end{array} \right] .
\]
Then $A=-A^{\top}$, and using $y_1 x_1+ \dots +y_n x_n =0$, it can be
verified that $Ax=y$. Finally, we note that $\|A\| \leq C \|y\|/|x_k|
\leq C(n, \delta) \|y\|$.
\end{proof}
\goodbreak

Let us note that whenever $M$ is a matrix in $\mC^{n\times n}$ such
that $M=-M^{\top}$, then $(Mx)^{\top}x=0$ for every $x\in \mC^n$.  We
will also need the following result on polynomial interpolation:

\begin{lemma}
\label{lemma_interpolation}
Given $a_k\in\mD$, and $m_k\geq 1$, $\alpha_k^{(m)}\in\mC$, $0\leq m
\leq m_{k}-1$, $k=1, \dots,N$, there exists a polynomial $p$ such that
\begin{equation}
\label{interpolation_constraints}
p^{(m)}(a_k)= \alpha_k^{(m)}, \quad k=1, \dots,N, \quad 0\leq m \leq
m_{k}-1.
\end{equation}
\end{lemma}
\begin{proof} For $k\in \{1, \dots, N\}$, let
\begin{eqnarray*}
p_{k}(z)&:=& (z-a_1)^{m_1} \dots (z-a_{k-1})^{m_{k-1}} (z-a_{k+1})^{m_{k+1}} \dots (z-a_N)^{m_{N}},
\\
q_{k}(z)&:=& 1/p_{k}(z).
\end{eqnarray*}
Define for $1\leq k \leq N,\; 0\leq m \leq m_k-1$:
\[
P_{k,m}(z):= p_{k}(z) \frac{(z-a_k)^{m}}{m!} \sum_{l=0}^{m_{k}-1-m}
\frac{q_{k}^{(l)}(a_k) }{l!} (z-a_k)^{l},
\]
Let
\[
p(z):= \sum_{k=1}^{N} \sum_{m=0}^{m_{k}-1} \alpha_{k}^{(m)} P_{k,m}(z).
\]
That $p$ satisfies \eqref{interpolation_constraints} can be verified
in the same manner as the proof of Theorem 1 in \cite{Spi60}.
\end{proof}

With these two lemmas, it is now possible to give the proof of
Theorem~\ref{theorem_estimates}.

Recall that for a vector $f=(f_1,\dots,f_n)\in (\H)^n$, we define its
norm by $||f||_\infty:=\max_{1\leq j\leq n}||f_j||_\infty$. Similarly
for a matrix $M=(f_{ij})_{1\leq i,j\leq n}\in (\H)^{n\times n}$.

\begin{proof}[Proof of Theorem \ref{theorem_estimates}]
By Carleson's Corona theorem, there exists a  $g_0 \in (H^\infty)^n$ such that
\[
g_0\cdot f =1,
\]
and $\|g_0\|_{\infty} \leq C( \delta)$. We will find a suitable
matrix $H \in (H^\infty)^{n\times n}$ such that
\begin{enumerate}
\item $H=-H^{\top}$,
\item (with norm control) $\|H\|_{\infty} \leq C(n, \delta, a_k,m_k)$, and
\item $g:= g_0 +Hf \in (H^\infty_B)^n $.
\end{enumerate}
Then clearly $g\cdot f=1$, giving Theorem \ref{theorem_estimates}.

Let $k\geq 1$. Since $g_0\cdot f =1$ and $f(a_k)=f(a_1)$,  we have
\[
g_0(a_1)\cdot f(a_1)=1 \;\;\textrm{ and }\;\; g_0(a_k)\cdot f(a_1)=1.
\]
Thus
\[
(g_0(a_k)-g_0(a_1))\cdot f(a_1)=0.
\]
Note that with $g:=g_0+Hf$, we have
\[
g(a_1)=g(a_k) \quad \textrm{iff} \quad
(H(a_k)-H(a_1))f(a_1) =g_0(a_1)-g_0(a_k).
\]
Choose $A_{a_1}^{(0)} \in \mC^{n \times n}$ such that
\[
(A_{a_1}^{(0)})^{\top}=-A_{a_1}^{(0)}.
\]
Next choose $A_{a_k}^{(0)} \in \mC^{n \times n}$ such that
\[
(A_{a_k}^{(0)})^{\top}=-A_{a_k}^{(0)}\;\; \textrm{ and }\;\;
(A_{a_k}^{(0)} -A_{a_1}^{(0)})f(a_1)=g_0(a_1)-g_0(a_k),
\]
which is possible in light of Lemma \ref{lemma_LA} and the choice of
$A_{a_1}^{(0)}$. We also note that since $\|g_0\|_{\infty} \leq
C(\delta)$, we can estimate $\|A_{a_k}^{(0)}\| \leq C(n, \delta)$.

Since $g_0\cdot f=1$, by differentiating this $m$ times, with $1\leq m\leq
m_k -1$, and evaluating at $a_k$, we obtain
\[
\sum_{j=0}^{m} \binom{m}{j} g_0^{(m-j)}(a_k)\cdot f^{(j)}(a_k)=0.
\]
Thus
\[
g_0^{(m)}(a_k)\cdot f(a_k)=0.
\]
We note that with $g:=g_0+Hf$, we have
\[
g^{(m)}(a_k)=0 \quad \textrm{if and only if} \quad
H^{(m)}(a_k)f(a_k)=-g_0^{(m)}(a_k).
\]
Again, using Lemma \ref{lemma_LA}, we can choose $A_{a_k}^{(m)} \in \mC^{n
  \times n}$ such that
\[
(A_{a_k}^{(m)})^{\top}=-A_{a_k}^{(m)}
\quad \textrm{and} \quad
A^{(m)}_{a_k}f(a_k)=-g_0^{(m)}(a_k).
\]
We note that by Cauchy integral formula, $|g_0^{(m)} (a_k)| \leq C(a_k,m_k)
\|g_0\|_{\infty}$, and so by Lemma \ref{lemma_LA},
$\|A^{(m)}_{a_k}\|\leq C(n, \delta,a_k,m_k)$.

We can find $\frac{n(n-1)}{2}$
  functions $h_{ij} \in H^\infty$, $2\leq i \leq j \leq n$ such that
\[
h_{ij}^{(m)}(a_k)=\big[A_{a_k}^{(m)}\big]_{ij}, \quad 0\leq m \leq
m_k -1 , \; k=1,\dots, n
\]
and
\[
\|h_{ij}\|_{\infty} \leq C(\sup \{\left|\big[A_{a_k}^{(m)}\big]_{ij}\right|\;:\;0\leq
m\leq m_{k}-1, \;k=1,\dots, n\}).
\]
Indeed, interpolating polynomials as in Lemma
\ref{lemma_interpolation} above will do.  Define
\[
H:= \left[ \begin{array}{ccccc}
0        & h_{1,2}  & h_{1,3} & \dots & h_{1,n} \\
-h_{1,2}  & 0       & h_{2,3} & \dots & h_{2,n} \\
-h_{1,3}  & -h_{2,3} & 0      & \dots & h_{3,n} \\
\vdots   & \vdots & \ddots & \ddots & \vdots \\
-h_{1, n-1} & \dots & \dots & 0      & h_{n-1,n} \\
-h_{1, n} & \dots & \dots & -h_{n-1, n} & 0
\end{array} \right] \in (H^\infty)^{n\times n}.
\]
Then $g:=g_0 +H f \in H^\infty_B$, $g\cdot f=1$ and
$\|g\|_{\infty} \leq C(n, \delta,a_k,m_k)$.
\end{proof}

\subsection{Generalisation to arbitrary ideals}

Recall that $H^\infty_{\mathbb{I}}$ is a subalgebra of $H^\infty$
given by $\mathbb{I}$, an ideal in $H^\infty$, with
\[
H^\infty_{\mathbb{I}}:=\{c+\varphi\;|\;
c\in \mC \textrm{ and } \varphi\in \mathbb{I}\}.
\]
We now will show that the Corona Theorem holds for these algebras.

\begin{proof}[Proof of Theorem \ref{theorem_estimates2}]
Let $f_{k}=c_{k}+\varphi_{k}$, where $c_{k} \in \mC$ and
$\varphi_{k} \in \mathbb{I}$ for all $k \in \{1, \dots, n\}$.  Since
$\mathbb{I}\not=\H$, we see that the ideal
$(\varphi_1,\cdots,\varphi_n)$ generated by the $\varphi_j$'s is
proper; hence $\inf_{z\in\mD}\sum_{j=1}^n |\varphi_j(z)|=0$. Thus
there exists a sequence $(z_k)_{k\in \mN}$ in $\mD$ so that
$\varphi_j(z_k)\to 0$ for every $j$.  Hence
$$
\delta\leq \sum_{j=1}^n|f_j(z_k)|=
\sum_{j=1}^n |c_j+\varphi_j(z_k)|\to \sum_{j=1}^n |c_j|.
$$
Therefore \eqref{CC_I} gives
\[
1\geq |c_1|+\dots+|c_n| \geq \delta.
\]
So all $c_k$'s satisfy $|c_k| \leq 1$ and at least one of the $c_k$'s
(which we may assume is $k=1$ without loss of generality) satisfies
$|c_k|\geq \delta/n$.  By Carleson's Corona Theorem, there exist
$\widetilde{g}_1, \dots, \widetilde{g}_n \in H^\infty$ such that
\[
1=f_1 \widetilde{g}_1 +\dots +f_n \widetilde{g}_n \textrm{ on } \mD,
\]
and $\|\widetilde{g}_k\|_{\infty} \leq C(\delta)$ for all $k$'s. Thus,
\begin{eqnarray*}
  1
  &=&
  \frac{1}{c_1} (f_1 -\varphi_1 ) = \frac{1}{c_1} (f_1 - 1\varphi_1 )
  =
  \frac{1}{c_1} \bigg(f_1 - \bigg(\sum_{k=1}^{n}f_{k}
  \widetilde{g}_k \bigg)\varphi_1\bigg)
  \\
  &=&
  f_1 \bigg(\frac{1}{c_1}-\widetilde{g}_1 \varphi_1 \bigg)
  + f_2 \bigg(-\frac{1}{c_1}\widetilde{g}_2 \varphi_1\bigg)
  + \dots+ f_n \bigg(-\frac{1}{c_1}\widetilde{g}_n \varphi_1\bigg)
  \\
  &=&f_1 g_1 +f_2 g_2 +\dots+f_n g_n,
\end{eqnarray*}
where
\[
g_1 :=\frac{1}{c_1}-\widetilde{g}_1 \varphi_1  ,\quad
g_2 := -\frac{1}{c_1}\widetilde{g}_2 \varphi_1 ,\quad
\dots, \quad
g_n := -\frac{1}{c_1}\widetilde{g}_n \varphi_1 .
\]
Clearly $g_1, \dots, g_n \in H^\infty_{\mathbb{I}}$, and since we have
$\frac{1}{c_1} \leq \frac{n}{\delta}$ and
\[
\|\varphi_1 \|_{\infty} =\|f_1-c_1 \|_{\infty} \leq
\|f_1\|_{\infty}+|c_1| \leq 1+1=2,
\]
it follows that $\|g_k\|_{\infty} \leq C(n,\delta)$ for all $k$.
\end{proof}

\section{Stable ranks}
\label{stableranks}

The notion of stable rank of a ring (which we call Bass stable rank)
was introduced by H.~Bass \cite{Bas64} to facilitate computations in
algebraic K-theory. We recall the definition of the Bass stable rank
of a ring below.

\begin{definition}
Let $A$ be an algebra with a unit $e$.  An n-tuple $a\in A^n$ is called \textit{unimodular} if there exists an n-tuple $b\in A^n$ such that $\sum_{j=1}^{n}a_jb_j=e$.  An n-tuple $a$ is called \textit{stable} or \textit{reducible} if there exists an (n-1)-tuple $x$ such that the (n-1)-tuple $(a_1+x_1a_n,\ldots, a_{n-1}+x_{n-1}a_{n})$ is unimodular.  The \textit{stable rank} (also called \textit{bsr}(A) in the literature) of the algebra $A$ is the least integer n such that every unimodular (n+1)-tuple is reducible.
\end{definition}


It is straightforward to see that for the algebras $A(\mD)$ and $H^\infty(\mD)$ that unimodular tuples are the same as Corona data in the algebra.

Over the last two decades there has been some interest in studying the
Bass stable rank of Banach algebras. Jones, Marshall and Wolff
\cite{JonMarWol} showed that the Bass stable rank of the disk algebra
$A(\mD)$ is one, and Treil \cite{Tre92} proved that the Bass stable
rank of $H^\infty$ is one as well; see also \cite{Sua94}:

\begin{proposition}[Jones-Marshall-Wolff]
The Bass stable rank of $A(\mD)$ is $1$.
\end{proposition}

\begin{proposition}[Treil]
The Bass stable rank of $H^\infty$ is $1$.
\end{proposition}

Using Treil's result, an analogous result holds for the
subalgebra $H^\infty_{\mathbb{I}}=\mC+\mathbb{I}$:

\begin{theorem}
  Let $\mathbb{I}$ be a proper ideal in $H^\infty$. Then the the Bass
  stable rank of $H^\infty_{\mathbb{I}}$ is $1$.
\end{theorem}
\begin{proof}
Let $(f,g)=(a+\varphi,b+\psi)$ be a unimodular pair in
$H^\infty_{\mathbb{I}}$, where $a,b\in\mC$ and $\varphi,\psi\in
\mathbb{I}$.  In particular $\delta:=\inf_{z\in\mD}|f(z)|+|g(z)|>0$.
We consider the two cases $a\neq 0$ and $a=0$.

\medskip

\noindent \underline{1}$^\circ$ If $a\neq 0$, then the pair,
$(f,g\varphi)$ is unimodular.  Indeed, if $z\in \mD$ is such that
$|\varphi(z)| \geq |a|/2$, then $|f|+|g\varphi| \geq \min \{1,
|a|/2\}(|f|+|g|) \geq \delta \min \{1, |a|/2\}$ at that $z$, while if
$|\varphi(z)|\leq |a|/2$, then $|f|=|a+\varphi| \geq |a|-|a|/2=|a|/2$,
so that $|f|+|g\varphi| \geq |a|/2$.  By Treil's Theorem, the element
$f+h(g\varphi)$ is invertible in $H^\infty$ for some element $h\in
H^\infty$. But writing $f+h(g\varphi)=f+(h\varphi)g$, and noting that
$h\varphi \in H^\infty_{\mathbb{I}}$, we are done since
$f+(h\varphi)g$ being invertible in $H^\infty$, is also invertible in
$H^\infty_{\mathbb{I}}$ by the corona theorem for
$H^\infty_{\mathbb{I}}$.

\medskip

\noindent \underline{2}$^\circ$ If $a=0$, then clearly $b\neq 0$, as
$(f,g)$ satisfies the corona condition. We know that there exist
$x,y\in H^\infty_{\mathbb{I}}$ such that $fx+gy=1$. Thus
$(f+g)x+g(y-x)=1$, and so $(f+g,g)$ is unimodular. By case 1$^\circ$
above, it follows that there exists a $h\in H^\infty_{\mathbb{I}}$
such that $(f+g)+hg$ is invertible in $H^\infty_{\mathbb{I}}$, that
is, $f+(1+h)g$ is invertible in $H^\infty_{\mathbb{I}}$. Since $1+h\in
H^\infty_{\mathbb{I}}$, we are done.
\end{proof}

By using the result of Jones, Marshall and Wolff that the Bass stable
rank of the disk algebra $A(\mD)$ is one, we now establish the following:

\begin{theorem}
  Let $B$ be a Blaschke product. Suppose that the zeros of $B$ do not
  cluster at a set of positive Lebesgue measure (hence $A(\mD)_B$ is not
  trivial). Then the Bass stable rank of $A(\mD)_B$ is $1$.
\end{theorem}

\begin{proof} Recall that $A(\mD)_B= (\mC +B\H)\cap A(\mD)$. Let $(f,g)=(a+BF,
  b+BG)$ be a unimodular pair in $A(\mD)_B$. Note that $F$ and $G$
  necessarily are in $A(\mD)$.  As above, we consider two cases:

\medskip

\noindent \underline{1}$^\circ$ If $a\neq 0$, then we let $E$ be the
set of cluster points of the zero set of $B$ (if $B$ is a finite Blaschke product,
then $E=\emptyset$).  Since $E$ has Lebesgue measure zero, there
exists by \cite[p. 81]{ho} a function $p\in A(\mD)$ that satisfies $p=1$ on
$E$ and $|p|<1$ on $\overline \mD\setminus E$. (If $E=\emptyset$, we
take $p\equiv 0$.)  Then $B(1-p)\in A(\mD)_B$.  Consider the pair $(f,
B(1-p)g)$. Then this is unimodular, too. Since $A(\mD)$ has the stable rank
$1$, there exists $h\in A(\mD)$ such that $f+h B(1-p)g $ is invertible in
$A(\mD)$.  But $B(1-p)h\in A(\mD)_B$.  So $f+h B(1-p)g\in A(\mD)_B$. Thus the pair
$(f,g)$ is reducible in $A(\mD)_B$.

\medskip

\noindent \underline{2}$^\circ$ If $a=0$, then we consider the
unimodular pair $(f+g,g)$ and use case 1$^\circ$ to conclude that
$(f+g,g)$ and hence $(f,g)$ is reducible in $A(\mD)_B$.
\end{proof}

We remark at this point that it is possible to prove a general
algebraic proposition along the same lines.  Namely, let $\mathcal{A}$
be a commutative unital, normed algebra and let $\mathbb{I}$ be a
non-zero ideal in $\mathcal{A}$.  If $\mathcal{A}$ has Bass stable
rank $n$, then $\mathbb{C}+\mathbb{I}$ has Bass stable rank $n$ as
well.

\section{Concluding Remarks}

Analysing the proof Theorem \ref{theorem_estimates2}, it seems
reasonable to make the following conjecture, where the new feature is
that the constant $C(\delta)$ is independent of the number $n$ of
generators. This conjecture is inspired by the fact that the
corresponding estimates in $H^\infty$ are true.

\begin{conjecture}
  Let $B$ be a Blaschke product with infinitely many zeros $a_k$ of
  multiplicity $m_k$. Suppose that $f_1,\ldots, f_n\in H^\infty_B$
  satisfy
\[
\forall z \in \mD, \;\; 1\geq\sum_{k=1}^{n} |f_k(z)| \geq \delta >0.
\]
Then there exist $g_1, \dots, g_n \in H^\infty_B$ such that
\[
\forall z\in \mD, \;\; \sum_{k=1}^n g_k (z)f_k(z)=1 \quad \textrm{and}
\quad \forall k \in \{1, \dots,
n\},\;\; \|g_{k}\|_{\infty} \leq C(\delta).
\]
\end{conjecture}

{\bf Acknowledgements.}  The authors wish to thank a careful referee.
The referee's comments helped in the overall presentation of the
paper. The third author thanks the University Paul Verlaine of Metz
for support during a one week research stay at the Laboratoire de
Math\'ematiques et Applications.



\end{document}